\documentclass[reqno,12pt,a4paper]{amsart}

\usepackage{mathrsfs}
\usepackage{amssymb}
\usepackage{amsfonts}
\usepackage{latexsym}
\usepackage{amsthm}

\usepackage{graphicx}
\def\lb{\label}

\newcommand{\er}[1]{\textrm{(\ref{#1})}}

\begin{document}


\renewcommand{\theequation}{\arabic{section}.\arabic{equation}}
\theoremstyle{plain}
\newtheorem{theorem}{\bf Theorem}[section]
\newtheorem{lemma}[theorem]{\bf Lemma}
\newtheorem{corollary}[theorem]{\bf Corollary}
\newtheorem{proposition}[theorem]{\bf Proposition}
\newtheorem{definition}[theorem]{\bf Definition}
\newtheorem{remark}[theorem]{\it Remark}

\def\a{\alpha}  \def\cA{{\mathcal A}}     \def\bA{{\bf A}}  \def\mA{{\mathscr A}}
\def\b{\beta}   \def\cB{{\mathcal B}}     \def\bB{{\bf B}}  \def\mB{{\mathscr B}}
\def\g{\gamma}  \def\cC{{\mathcal C}}     \def\bC{{\bf C}}  \def\mC{{\mathscr C}}
\def\G{\Gamma}  \def\cD{{\mathcal D}}     \def\bD{{\bf D}}  \def\mD{{\mathscr D}}
\def\d{\delta}  \def\cE{{\mathcal E}}     \def\bE{{\bf E}}  \def\mE{{\mathscr E}}
\def\D{\Delta}  \def\cF{{\mathcal F}}     \def\bF{{\bf F}}  \def\mF{{\mathscr F}}
\def\c{\chi}    \def\cG{{\mathcal G}}     \def\bG{{\bf G}}  \def\mG{{\mathscr G}}
\def\z{\zeta}   \def\cH{{\mathcal H}}     \def\bH{{\bf H}}  \def\mH{{\mathscr H}}
\def\e{\eta}    \def\cI{{\mathcal I}}     \def\bI{{\bf I}}  \def\mI{{\mathscr I}}
\def\p{\psi}    \def\cJ{{\mathcal J}}     \def\bJ{{\bf J}}  \def\mJ{{\mathscr J}}
\def\vT{\Theta} \def\cK{{\mathcal K}}     \def\bK{{\bf K}}  \def\mK{{\mathscr K}}
\def\k{\kappa}  \def\cL{{\mathcal L}}     \def\bL{{\bf L}}  \def\mL{{\mathscr L}}
\def\l{\lambda} \def\cM{{\mathcal M}}     \def\bM{{\bf M}}  \def\mM{{\mathscr M}}
\def\L{\Lambda} \def\cN{{\mathcal N}}     \def\bN{{\bf N}}  \def\mN{{\mathscr N}}
\def\m{\mu}     \def\cO{{\mathcal O}}     \def\bO{{\bf O}}  \def\mO{{\mathscr O}}
\def\n{\nu}     \def\cP{{\mathcal P}}     \def\bP{{\bf P}}  \def\mP{{\mathscr P}}
\def\r{\rho}    \def\cQ{{\mathcal Q}}     \def\bQ{{\bf Q}}  \def\mQ{{\mathscr Q}}
\def\s{\sigma}  \def\cR{{\mathcal R}}     \def\bR{{\bf R}}  \def\mR{{\mathscr R}}
\def\S{\Sigma}  \def\cS{{\mathcal S}}     \def\bS{{\bf S}}  \def\mS{{\mathscr S}}
\def\t{\tau}    \def\cT{{\mathcal T}}     \def\bT{{\bf T}}  \def\mT{{\mathscr T}}
\def\f{\phi}    \def\cU{{\mathcal U}}     \def\bU{{\bf U}}  \def\mU{{\mathscr U}}
\def\F{\Phi}    \def\cV{{\mathcal V}}     \def\bV{{\bf V}}  \def\mV{{\mathscr V}}
\def\P{\Psi}    \def\cW{{\mathcal W}}     \def\bW{{\bf W}}  \def\mW{{\mathscr W}}
\def\o{\omega}  \def\cX{{\mathcal X}}     \def\bX{{\bf X}}  \def\mX{{\mathscr X}}
\def\x{\xi}     \def\cY{{\mathcal Y}}     \def\bY{{\bf Y}}  \def\mY{{\mathscr Y}}
\def\X{\Xi}     \def\cZ{{\mathcal Z}}     \def\bZ{{\bf Z}}  \def\mZ{{\mathscr Z}}
\def\O{\Omega}

\newcommand{\gA}{\mathfrak{A}}
\newcommand{\gB}{\mathfrak{B}}
\newcommand{\gC}{\mathfrak{C}}
\newcommand{\gD}{\mathfrak{D}}
\newcommand{\gE}{\mathfrak{E}}
\newcommand{\gF}{\mathfrak{F}}
\newcommand{\gG}{\mathfrak{G}}
\newcommand{\gH}{\mathfrak{H}}
\newcommand{\gI}{\mathfrak{I}}
\newcommand{\gJ}{\mathfrak{J}}
\newcommand{\gK}{\mathfrak{K}}
\newcommand{\gL}{\mathfrak{L}}
\newcommand{\gM}{\mathfrak{M}}
\newcommand{\gN}{\mathfrak{N}}
\newcommand{\gO}{\mathfrak{O}}
\newcommand{\gP}{\mathfrak{P}}
\newcommand{\gQ}{\mathfrak{Q}}
\newcommand{\gR}{\mathfrak{R}}
\newcommand{\gS}{\mathfrak{S}}
\newcommand{\gT}{\mathfrak{T}}
\newcommand{\gU}{\mathfrak{U}}
\newcommand{\gV}{\mathfrak{V}}
\newcommand{\gW}{\mathfrak{W}}
\newcommand{\gX}{\mathfrak{X}}
\newcommand{\gY}{\mathfrak{Y}}
\newcommand{\gZ}{\mathfrak{Z}}

\def\ve{\varepsilon}   \def\vt{\vartheta}    \def\vp{\varphi}    \def\vk{\varkappa}

\def\Z{{\mathbb Z}}    \def\R{{\mathbb R}}   \def\C{{\mathbb C}}    \def\K{{\mathbb K}}
\def\T{{\mathbb T}}    \def\N{{\mathbb N}}   \def\dD{{\mathbb D}}


\def\la{\leftarrow}              \def\ra{\rightarrow}            \def\Ra{\Rightarrow}
\def\ua{\uparrow}                \def\da{\downarrow}
\def\lra{\leftrightarrow}        \def\Lra{\Leftrightarrow}


\def\lt{\biggl}                  \def\rt{\biggr}
\def\ol{\overline}               \def\wt{\widetilde}
\def\no{\noindent}


\let\ge\geqslant                 \let\le\leqslant
\def\lan{\langle}                \def\ran{\rangle}
\def\/{\over}                    \def\iy{\infty}
\def\sm{\setminus}               \def\es{\emptyset}
\def\ss{\subset}                 \def\ts{\times}
\def\pa{\partial}                \def\os{\oplus}
\def\om{\ominus}                 \def\ev{\equiv}
\def\iint{\int\!\!\!\int}        \def\iintt{\mathop{\int\!\!\int\!\!\dots\!\!\int}\limits}
\def\el2{\ell^{\,2}}             \def\1{1\!\!1}
\def\sh{\sharp}
\def\wh{\widehat}
\def\bs{\backslash}

\def\where{\mathop{\mathrm{where}}\nolimits}
\def\all{\mathop{\mathrm{all}}\nolimits}
\def\as{\mathop{\mathrm{as}}\nolimits}
\def\Area{\mathop{\mathrm{Area}}\nolimits}
\def\arg{\mathop{\mathrm{arg}}\nolimits}
\def\const{\mathop{\mathrm{const}}\nolimits}
\def\det{\mathop{\mathrm{det}}\nolimits}
\def\diag{\mathop{\mathrm{diag}}\nolimits}
\def\diam{\mathop{\mathrm{diam}}\nolimits}
\def\dim{\mathop{\mathrm{dim}}\nolimits}
\def\dist{\mathop{\mathrm{dist}}\nolimits}
\def\Im{\mathop{\mathrm{Im}}\nolimits}
\def\Iso{\mathop{\mathrm{Iso}}\nolimits}
\def\Ker{\mathop{\mathrm{Ker}}\nolimits}
\def\Lip{\mathop{\mathrm{Lip}}\nolimits}
\def\rank{\mathop{\mathrm{rank}}\limits}
\def\Ran{\mathop{\mathrm{Ran}}\nolimits}
\def\Re{\mathop{\mathrm{Re}}\nolimits}
\def\Res{\mathop{\mathrm{Res}}\nolimits}
\def\res{\mathop{\mathrm{res}}\limits}
\def\sign{\mathop{\mathrm{sign}}\nolimits}
\def\span{\mathop{\mathrm{span}}\nolimits}
\def\supp{\mathop{\mathrm{supp}}\nolimits}
\def\Tr{\mathop{\mathrm{Tr}}\nolimits}
\def\BBox{\hspace{1mm}\vrule height6pt width5.5pt depth0pt \hspace{6pt}}
\def\ch{\mathop{\mathrm{ch}}\nolimits}


\newcommand\nh[2]{\widehat{#1}\vphantom{#1}^{(#2)}}
\def\dia{\diamond}

\def\Oplus{\bigoplus\nolimits}



\def\qqq{\qquad}
\def\qq{\quad}
\let\ge\geqslant
\let\le\leqslant
\let\geq\geqslant
\let\leq\leqslant
\newcommand{\ca}{\begin{cases}}
\newcommand{\ac}{\end{cases}}
\newcommand{\ma}{\begin{pmatrix}}
\newcommand{\am}{\end{pmatrix}}
\renewcommand{\[}{\begin{equation}}
\renewcommand{\]}{\end{equation}}
\def\eq{\begin{equation}}
\def\qe{\end{equation}}
\def\[{\begin{equation}}
\def\bu{\bullet}


\newcommand{\ee}{\epsilon}
\newcommand{\abs}[1]{\lvert#1\rvert}
\newcommand{\norm}[1]{\lVert#1\rVert}
\renewcommand{\SS}{{\mathfrak S}}

\title[{Sharp asymptotics of the quasimomentum}]
{Sharp asymptotics of the quasimomentum}

\date{\today}
\author[Evgeny L. Korotyaev]{Evgeny L. Korotyaev}
\address{
Saint-Petersburg University; Saint-Petersburg, Russia,
 \ korotyaev@gmail.com}

\date{\today}


\begin{abstract}

We consider the Schr\"odinger operator with a periodic potential $p$
on the real line. We assume that $p$ belongs to the Sobolev space
$\mH_m$ on the circle for some  $m\ge -1$, and we determine the
asymptotics of the quasimomentum and the Titchmarsh-Weyl functions,
the Bloch functions at high energy.

\end{abstract}

\maketitle

\section {Introduction and main results}

\setcounter{equation}{0}

Consider the Schr\"odinger operator $H$ acting in the Hilbert space
$L^2(\R)$  and  given by
$$
 Hf=-f''+pf.
$$
Here the potential $p$ is 1-periodic and belongs to  the Sobolev
space $\mH_m$ on the circle $\T={\R/\Z}$:
\[
\lb{1}
p\in \mH_m=\{p^{(m)}\in L^2(\T)\}, \qqq m\ge -1.
\]
We recall the results from \cite{K1} about the operator $H$. The
spectrum of $H$ is absolutely continuous and has the form $\s(H_0)=
\bigcup\limits_{n\in\N} \gS_n$,
 where the bands $\gS_n$ and gaps $\g_n$ are given by
$$
\gS_n=[E^+_{n-1},E^-_n],\ \ \qq \g_{n}=(E^-_{n},E^+_n), \qq \forall
n\in\N=\{n: n=1,2,3,....\},
$$
see Fig. 1.
Without loss of generality, we may assume $E_0^+=0$. Here the $E_n^\pm$ satisfy
\[
\lb{2}
0=E_0^+< E^-_1 \le E^+_1\dots\le E^+_{n-1}< E^-_n \le E^+_{n}<\dots
\]
If $p\in \mH_m$, then it is known that there are infinitely many
non-degenerate gaps, i.e. $E_n^- < E_n^+$, unless $p$ is arbitrarily
often differentiable, and all gaps are non-degenerate generically
(see e.g. \cite{MO}, \cite{K1}). The sequence \er{2} is the spectrum
of the equation
\[
\lb{3} -y''+py=\l y,
\]
with the condition of 2-periodicity, $y(x+2)=y(x)$ $(x\in \R)$. If a
gap degenerates, $\g_n=\es $ for some $n$, then the corresponding
bands $\gS_{n} $ and $\gS_{n+1}$ touch. This happens when
$E_n^-=E_n^+$; this number is then a double eigenvalue of the
2-periodic problem \er{3}. The lowest eigenvalue $E_0^+=0$ is always
simple and has a 1-periodic eigenfunction. Generally, the
eigenfunctions corresponding to the eigenvalues $E_{2n}^{\pm}$ are
1-periodic, and those for $E_{2n+1}^{\pm}$ are 1-anti-periodic in
the sense that $y(x+1)=-y(x)$ $(x\in\R)$.

In the case of the potential $p\in \mH_m, m\ge 0$, throughout the
paper,  we shall denote by $\vt(x,z)$, $\vp(x,z)$ the two solutions
forming the canonical fundamental system of the unperturbed equation
\[
-y''+py=z^2y,
\]
under the initial conditions
$$
\vp'(0,z)=\vt(0,z)=1, \qqq  \qqq \vp(0,z)=\vt'(0,z)=0.
 $$
Here and in the following $"\ ' \ "$ denotes the derivative w.r.t.
the first variable. In the following, we shall treat the momentum
$z=\sqrt \l$ (as opposed to the energy $z^2=\l$) as the principal
spectral variable. The Lyapunov function (which is the Hill
discriminant for $m\ge 0$) of the periodic equation is then defined
by
$$
\D(z)={1\/2}(\vp'(1,z)+\vt(1,z)).
$$
In the case $m=-1$ we denote the  Lyapunov function also by $\D(z)$.
In the last case the definition of $\D(z)$ is more complicated and
is given in Section 4. Recall that the function $\D(z)$ is entire
and even $\D(-z)=\D(z), z\in \cZ$.

We introduce the {\it quasimomentum} $k(\cdot)$ for $H$ as
$k(z)=\arccos \D(z), z \in \cZ$, where $\cZ$ is the cut domain (see
Fig. 1 and 2) given by
\[
\lb{5} \cZ=\C\sm \bigcup_{n\in \Z} \ol g_n,\qq {\rm where} \qq
g_n=(e_n^-,e_n^+)=-g_{-n},\qq e_n^\pm=\sqrt{E_n^\pm}>0,\qq n\ge 1,
\qq g_0=\es.
\]
Note that $\D(e_{n}^{\pm})=(-1)^n$ and if $\l\in \g_n,  n\ge 1$,
then $z\in g_{\pm n}$, and if $\l\in \g_0=(-\iy,E_0^+)$, then $z\in
i\R_\pm$. The function $k$ is analytic in $\cZ$ and satisfies
\[
\begin{aligned}
\lb{pk}
&(i)\qqq k(z)=z+o(1)\qqq  as \ \ \Im z\to \iy,
\qqq \qqq \qqq \qqq \qqq \qqq \qqq \\
& (ii)\qqq   k(0)=0,\qqq
\Re k(z\pm i0)|_{[e_n^-,e_n^+]}=\pi n\qq (n\in \Z),
\qqq \qqq \qqq \qqq \qqq  \\
& (iii)\qqq k(-z)=-k(z), \qqq  \forall \qqq z\in \cZ,\\
& (iv)\qqq \pm \Im k(z)>0, \qqq\forall \ z\in \C_\pm=\{z\in \C: \pm \Im z>0  \},
\end{aligned}
\]
see (\cite{MO}, \cite{KK}). Moreover, $k$ is a conformal mapping
from $\cZ$ onto the  quasimomentum domain $\cK$ given by
\[
\lb{cK}
\cK=\C\sm \cup \ol\G_n, \qqq \G_n=(\pi n-ih_n,\pi n+ih_n),
\]
see  Figs. 2 and 3. Here $\G_n$ is a vertical cut of the height
$h_n=h_{-n}\ge 0, h_0=0$. The height $h_n$ is determined by the
equation $\cosh h_n=|\D(e_n)|\ge 1$,  where $e_n\in [e_n^-,e_n^+]$
is such that $\D'(e_n)=0$. Note that the point $e_n$ is  unique for
each $n\in \Z$. The function $k$ maps the cut $g_n$ onto the cut
$\G_n$.

We have obtained a conformal mapping $k: \cZ\to \cK$, called {\it
the quasimomentum mapping} (or shortly the quasimomentum), which
generalizes the classical quasimomentum (see e.g. \cite{RS}). A
point $z\in \cZ$ is called {\it a momentum} and a point $k\in \cK$
is called {\it a quasimomentum}. The abstract quasimomentum, which
we have just defined is related to the spectral theory of the Hill
operator $H$ by the following construction invented in \cite{F1},
\cite{F}, \cite{MO} for the $L^2$ potentials and generalized in
\cite{K1} for the potential from $\mH_{-1}$. Some asymptotics of the
quasimomentum for $p\in \mH_0$ were obtained in \cite{F2}, outside
some  neighborhoods of gaps. The quasimomentum for the Schr\"odinger
operator $-{d^2\/dx^2}+V$ acting on the real line where $V$ is a
periodic $N\ts N$ matrix-valued potential was studied in \cite{CK}.
We would like to add that the properties of the quasimomentum are
important in many different fields, see e.g. : inverse problem
\cite{F},   \cite{GT}, \cite{KK1},  \cite{K1}, \cite{MO}, non-linear
equations \cite{C}, \cite{GWH}, and so on.

For any  $p\in \mH_m, m\ge 0$ we define the integrals
\[
\begin{aligned}
\lb{Pj}
P_{-1}={\int_0^1pdx\/2},\qq  P_0={\int_0^1p^2dx\/2^3},\qq
P_{j}={\|p^{(j)}\|^2+\int_0^1F_jdx\/2^{3+2j}},\qq j=1,...,m.
\end{aligned}
\]
Here $F_j$ is some polynomial of $p,p',p'',\dots ,p^{(j-1)}$.
 In particular, we have
\[
\begin{aligned}
\lb{6.6}
\qqq F_1=2p^3, \qqq F_2=10p{p'}^2+5p^4,\qqq
F_3=14p{p''}^2+70p^2{p'}^2+112p^5,\dots ,
\end{aligned}
\]
see \cite{MM}, \cite{MO}, where all $P_j>0$ if $p\ne 0$ and
$E_0^+=0$, since we have \er{QmP}. Introduce the functions
\[
\lb{Km}
 K_m(z)={P_{-1}\/z}+{P_{0}\/z^{3}}+...+{P_{m-1}\/z^{2m+1}},
\]
and define the domains
$$
\cZ_\ve=\{z\in \cZ:  \dist\{z,g\}>\ve\}, \qq \ve>0,\qqq \where \qq
g=\bigcup_{n\in \Z}g_n.
$$

\medskip

\begin{theorem}
\lb{T1} Let $p\in \mH_{m}$ for
some $m\ge 0$ and let $A, \ve>0$.  Then
\[
\lb{13}
k=z-K_{m}(z)+f_{m+1}(z),\qq f_{m+1}(z)={1\/\pi z^{2m+2}}
\int_{\R}{t^{2m+2}v(t)dt\/ t-z},\qqq \ z\in \cZ,
\]
where $f_{m+1}$ has  the following asymptotics as $|z|\to \iy$:
\[
\lb{ask1}
f_{m+1}(z)=-{P_{m}+o(1)\/z^{2m+3}}\qqq as \qqq z\in \{z=x+iy\in\C: y>A|x|\},
\]
\[
\lb{ask2}
f_{m+1}(z)={O(1)\/z^{2m+2}} \qq \as \qq z\in \cZ_\ve,
\]
\[
\lb{ask3}
|f_{m+1}(z)|\le {|\g_n|\/2\pi n}+b_n,\qqq
b_n={O(|\g_n|)\/n^3}+{O(1)\/n^{2m+2}}\qqq \dist\{z, g_n\}\le \ve.
\]
Moreover, the asymptotic estimate \er{ask3} is sharp, since
\[
\lb{ask4}
f_{m+1}(e_n^\pm)=\mp {|\g_n|\/2\pi n}(1+o(1)) \qq \as \qq n\to \iy.
\]
\end{theorem}

{\bf Remarks.} 1)  Recall that $p\in \mH_m$ if and only
if $(n^m|\g_n|)_1^\iy\in \ell^2$, see \cite{MO}, \cite{K3}.

2)  \er{ask1}-\er{ask3} give 3 types of  asymptotics. The "best"
asymptotics \er{ask1} has the form $f_{m+1}(z)={O(1)\/z^{2m+3}}$ and
the "bad"  asymptotics \er{ask3}  has the form
$f_{m+1}(z)={O(n^m|\g_n|)\/n^{m+1}}$. There is a big difference
between the sharp asymptotics \er{ask1} and \er{ask3}, since  due to
\er{ask4} the asymptotics \er{ask3}  is sharp.

 3)  Shenk and Shubin \cite{SS} determined complete asymptotic expansions
 of the integrated density of  states for $p\in C^\iy(\R)$. Recall
 that the integrated  density of  states is given by
 ${1\/\pi} \Re k(z), z\in \R$.

4) The asymptotics \er{ask1}-\er{ask3} give an  asymptotics of the
integrated density of   states for the case $p\in \mH_m$. If $p\in
C^\iy(\R)$, then the theorem gives  complete asymptotic expansions
of  the quasimomentum $k(z)$.

5) The complete asymptotic expansion of the integrated density of
 states of multidimensional almost-periodic Schrodinger operators
 were determined by Parnovski, Shterenberg \cite{PS1},
 see also \cite{KP}, \cite{PS2} and  the references therein.

\medskip

In Section 4 we consider the case of distributional potentials
$p\in\mH_{-1}$.

\medskip

In order to write the more complete results about the asymptotics
for the Hill operator, we determine the asymptotics of the Bloch
functions and the Titchmarsh-Weyl function. Note that, although
asymptotic expressions for the Bloch functions  and the
Titchmarsh-Weyl function for the case $p\in \mH_m$ have not been
formally written out anywhere   previously,  this  result can be
regarded as known, since it can easily be obtained with the help of
the results of Marchenko and Ostrovskii \cite{MO}.

We introduce the Bloch functions $\P_{\pm}$ of $H$ defined by (see \cite{T})
\[
\P_{\pm}(x,z)=\vp(x,z)+M_\pm(z)\vt(x,z), \qqq (x,z)\in [0,1]\ts \cZ,
\]
where $M_{\pm}(z) $ is the Titchmarsh-Weyl function given by
\[
  M_{\pm}(z)={\b (z )\pm \sin k(z )\/ \vp(1,z )},\ \ \ \
\ \ \ \ \ \b(z )={\vp'(1,z )-\vt (1,z) \/2}.
\]
Furthermore, we introduce the model function (see Lemma 3.1
in\cite{MO})
\[
\begin{aligned}
\lb{ym}
 \x_m(x,z)=z x-i\int_0^x\sum_1^{m}{\vk_j(t)\/(2iz)^j}dt,
 \qqq (x,z)\in [0,1]\ts \C, \ \ z\ne 0.
\end{aligned}
\]
Here the functions $\vk_{j}$ are constructed with the help of
the recursion relations:
\[
\lb{ym1}
 \vk_{j+1}=-\vk_{j}'-\sum_1^{j-1} \vk_{j-s}\vk_s,\ \ j=1,2,...,m-1,
\]
where, in particular,
\[
\begin{aligned}
\lb{ym2}
\vk_1=p,\ \ \vk_2=-p',\ \ \
\vk_3=p''-p^2,\ \ \ \vk_4=-p'''+4pp',\ \ \ .... \\
\vk_j=(-1)^{j-1}p^{(j-1)}+\cP_{j-3}, \qq j=2,3,..., m,
\end{aligned}
\]
and $\cP_{j}$ is a polynomial in $p, p', p'', ..., p^{(j)}$.

\begin{theorem}
\lb{T2} Let $p\in \mH_{m}$ for some $m\ge 0$ and let $\ve>0, r\ge
1$. Assume that $E_0^+$ is any real number. Then the following
asymptotics hold true as $|z|\to \iy$:
\[
\lb{aM}
M_{\pm}(z)=i\x_m'(0,\pm z)+O(z^{1-m}),
\]
\[
\lb{aP}
\P_{\pm}(x,z)=e^{i \x_m(x,\pm z)}+O(z^{-m}),
\]
as $z\in \cZ_\ve, |\Im z|<r$, uniformly in $x\in [0,1]$.

Moreover,  if in addition $E_0^+=0$, then
\[
\begin{aligned}
\lb{aqu}
k(z)=\x(1,z)+O(z^{-m}),\ \
\\
{(-1)^{j}\/2^{2j+1}}\int_0^1\vk_{2j+1}(t)dt=Q_{2j},\ \ \ \ \
\int_0^1\vk_{2j}(t)dt=0,
\ \ \  j\ge 0,
\end{aligned}
\]
as $|\Im z|\le r$ and $|z|\to \iy$.

\end{theorem}

\no {\bf Remarks} 1)  Shenk and Shubin \cite{SS} determined
complete asymptotic expansions of the Bloch functions for
$p\in C^\iy(\R)$. There are some asymptotics of the Bloch
functions for $p\in L^1(0,1)$ in \cite{T}, \cite{F}, \cite{F2}.

2) In the proof of the theorem we use the standard asymptotics
of the solutions of
the equation $-y''+p(x)y=z^2 y$ for large $z$ from \cite{MO}.

\

\section {Asymptotics of the quasimomentum}
\setcounter{equation}{0}

\

Recall that  the quasimomentum $k(z)$ is a conformal mapping from
the momentum domain $\cZ$
 onto  the quasi-momentum domain $\cK$ given by (see Fig. 2 and 3)
\[
\begin{aligned}
\lb{mD}
\cZ=\C\sm \cup \ol g_n,\qq{\rm where} \qq
g_n=(e_n^-,e_n^+)=-g_{-n},\qq e_n^\pm=\sqrt{E_n^\pm}>0,\qq
n\ge 1,\qq g_0=\es,
\\
\cK=\C\sm\cup\ol\G_n, \qq {\rm where} \qq
\G_n=(\pi n+ih_n,\pi n-ih_n),\qq h_n=h_{-n}\ge 0,  \ \ n\ge 1, \ \ h_0=0.
\end{aligned}
\]
The height $h_n$ is determined by the equation $\cosh
h_n=|\D(e_n)|\ge 1$,  where $e_n\in [e_n^-,e_n^+]$ is such that
$\D'(e_n)=0$. Note that $e_n$ is  unique for each $n\in \Z$.
Cutting the $n$-th momentum gap  $g_n$ (if non-empty), we obtain a
cut $g_n^c$ with upper rim $g_n^+$ and lower rim $g_n^-$. Below, we
will identify this cut $g_n^c$ and the union of the upper rim (gap)
$\ol g_{n}^+$ and the lower rim (gap) $\ol g_{n}^{\ -}$,
 i.e.,
\[
g_n^c=\ol g_{n}^+\cup \ol g_{n}^-,\qqq {\rm where}\ g_{n}^\pm =g_n\pm i0;
\qq {\rm and}\ z\in g_n \Rightarrow z\pm i0\in g_n^\pm.
\]
Any non-degenerate (degenerate) cut $\G_n$ is connected in the some
way with the non-degenerate (degenerate) gap $\g_n$ and the
momentum gap $g_n$. We introduce the decomposition $k=u+iv$, where
$u,v$ are real harmonic functions in $\cZ$. The function $u(z)=\Re
k(z) $ is strongly  increasing on each band $\s_n$
 and equals $\pi n $ on each gap $[z_n^-,z_n^+],\ n\in\Z $;
the function $v(z)=\Im k(z) $ equals zero  on each  band $\s_n$,
is  strongly concave on each gap $g_n$ and has the maximum $h_n$ in $g_n$,
 attained  at some point
$e_n$, so that $h_n=v(e_n)$. Here and below we write
\[
\lb{vi0}
v(z)=v(z+i0) \qqq \as \qq  z\in \R.
\]
If $h_n=0$, then  n-the gap  is empty and  $e_n^-=e_n^+=e_n$.
These and  others properties of the comb mappings
can be found in \cite{KK},\cite{MO}.

Introduce the real  spaces
$$
 \ell^a=\rt\{f=(f_n)_{n\ge 1},\ \ \| f \|_a<\iy \rt\},
\qqq \ \ \| f \|_a^a=\sum _{n\ge 1} |f_n|^a <\iy, \ a\ge 1.
$$

Now we briefly discuss the properties of the general quasimomentum mapping
$k=u+iv$, as a function of $z=x+iy\in \cZ$. Their proof may be found in
 \cite{MO}, \cite{K1}, \cite{K2}, \cite{K4}.

\no 1) {\it $v(z)\ge \Im z>0$ and $v(z)=-v(\ol z)$ for all $z\in \C_+$ and
\[
\lb{2.3} k(-z)=-k(z)=\ol k(\ol z) ,\qq \all \ z\in \cZ.
\]
\no 2) $v(z)=0$ for all $z\in \s_n=[e_{n-1}^+,e_n^-], n\in \Z$.

\no 3) If some $g_n\ne \es, n\in \Z$, then $v(z)>0$ and $v''(z)<0$
for all $z\in g_n$, and $v(z)$ has a maximum at $e_n\in g_n$ such
that $v'(e_n)=0$, see Fig. 3,  and $\D'(e_n)=0$  and
\[
\lb{prqx}
v(z+i0)=-v(z-i0)>0, \qqq \qqq  \all \  z\in g_n\ne \es,
\]
\[
\lb{em-1}
|g_n|\le 2h_n,\qqq v(e_n)=h_n>0.
\]
 Recall that $v(z)=v(z+i0)$ for all $z\in \R.$

\no 4) $u'(z)>0$ on  all $(e_{n-1}^+,e_n^-)$ and  $u(z)=\pi n$
for all $z\in g_n\ne \es, n\in \Z$.

\no 5) The function $k(z)$ maps a horizontal cut (a "gap" )
$[e_n^-,e_n^+]$  onto the vertical cut $\ol\G_n$  and
a spectral band $\s_n$ onto the segment $[\pi (n-1), \pi n]$ for all $n\in\Z$.

\no 6) The following asymptotics hold true:
\[
\lb{aen}
e_n^\pm=\pi n+o(1) \qqq \as \qq n\to \iy.
\]
\no  7) The following identity holds true:
\[
\lb{Kq}
k(z)=z+{1\/\pi}\int_g {v(t)\/t-z}dt,\qqq \forall z\in \cZ, \qq g=\cup g_n.
\]
\no  8) Introduce the moments
$$
 Q_{m}={1\/\pi}\int_\R t^{m}v(t+i0)dt<\iy,\qqq m\ge 0
$$
and note that $ Q_{m}=0$ for odd $m\ge 1$. Then the  following
identities and estimate hold true:
\[
\lb{QmP}
 Q_{2m+2}=P_m,\qqq   m\ge -1,
\]
\[
\lb{hQ}
\|h\|_\iy^2\le 2Q_0.
\]
}

\bigskip

If $p\in \mH_0$, then the quasimomentum $k(\cdot)$ has
the asymptotics (see \cite{K2})
\[
\lb{ak}
k(z)=z-{Q_0\/z}-{Q_2+o(1)\/z^3}\qqq as \qq \Im z\to \iy.
\]

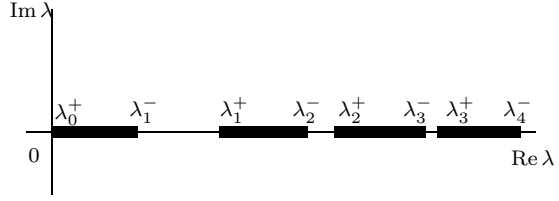
\begin{figure}
\tiny
\unitlength=1.00mm
\special{em:linewidth 0.4pt}
\linethickness{0.4pt}
\begin{picture}(108.67,33.67)
\put(41.00,17.33){\line(1,0){67.67}}
\put(44.33,9.00){\line(0,1){24.67}}
\put(108.33,14.00){\makebox(0,0)[cc]{$\Re\l$}}
\put(41.66,33.67){\makebox(0,0)[cc]{$\Im\l$}}
\put(42.00,14.33){\makebox(0,0)[cc]{$0$}}
\put(44.33,17.33){\linethickness{4.0pt}\line(1,0){11.33}}
\put(66.66,17.33){\linethickness{4.0pt}\line(1,0){11.67}}
\put(82.00,17.33){\linethickness{4.0pt}\line(1,0){12.00}}
\put(95.66,17.33){\linethickness{4.0pt}\line(1,0){11.00}}
\put(46.66,20.00){\makebox(0,0)[cc]{$\l_0^+$}}
\put(56.66,20.33){\makebox(0,0)[cc]{$\l_1^-$}}
\put(68.66,20.33){\makebox(0,0)[cc]{$\l_1^+$}}
\put(78.33,20.33){\makebox(0,0)[cc]{$\l_2^-$}}
\put(84.33,20.33){\makebox(0,0)[cc]{$\l_2^+$}}
\put(93.00,20.33){\makebox(0,0)[cc]{$\l_3^-$}}
\put(98.66,20.33){\makebox(0,0)[cc]{$\l_3^+$}}
\put(106.33,20.33){\makebox(0,0)[cc]{$\l_4^-$}}
\end{picture}
\caption{The spectral domain $\C\sm \cup \gS_n$ and the bands
$\gS_n=[\l^+_{n-1},\l^-_n], n\ge 1$}
\lb{sS}
\end{figure}
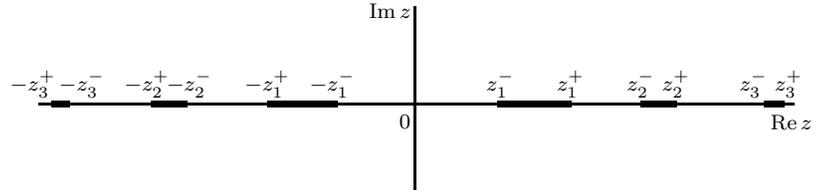
\begin{figure}
\tiny
\unitlength=1mm
\special{em:linewidth 0.4pt}
\linethickness{0.4pt}
\begin{picture}(120.67,34.33)
\put(20.33,21.33){\line(1,0){100.33}}
\put(70.33,10.00){\line(0,1){24.33}}
\put(69.00,19.00){\makebox(0,0)[cc]{$0$}}
\put(120.33,19.00){\makebox(0,0)[cc]{$\Re z$}}
\put(67.00,33.67){\makebox(0,0)[cc]{$\Im z$}}
\put(81.33,21.33){\linethickness{2.0pt}\line(1,0){9.67}}
\put(100.33,21.33){\linethickness{2.0pt}\line(1,0){4.67}}
\put(116.67,21.33){\linethickness{2.0pt}\line(1,0){2.67}}
\put(60.00,21.33){\linethickness{2.0pt}\line(-1,0){9.33}}
\put(40.00,21.33){\linethickness{2.0pt}\line(-1,0){4.67}}
\put(24.33,21.33){\linethickness{2.0pt}\line(-1,0){2.33}}
\put(81.67,24.00){\makebox(0,0)[cc]{$z_1^-$}}
\put(91.00,24.00){\makebox(0,0)[cc]{$z_1^+$}}
\put(100.33,24.00){\makebox(0,0)[cc]{$z_2^-$}}
\put(105.00,24.00){\makebox(0,0)[cc]{$z_2^+$}}
\put(115.33,24.00){\makebox(0,0)[cc]{$z_3^-$}}
\put(120.00,24.00){\makebox(0,0)[cc]{$z_3^+$}}
\put(59.33,24.00){\makebox(0,0)[cc]{$-z_1^-$}}
\put(50.67,24.00){\makebox(0,0)[cc]{$-z_1^+$}}
\put(40.33,24.00){\makebox(0,0)[cc]{$-z_2^-$}}
\put(34.67,24.00){\makebox(0,0)[cc]{$-z_2^+$}}
\put(26.00,24.00){\makebox(0,0)[cc]{$-z_3^-$}}
\put(19.50,24.00){\makebox(0,0)[cc]{$-z_3^+$}}
\end{picture}
\caption{$z$-domain $\cZ=\C\sm\cup g_n$, where $z=\sqrt{\l}$
and momentum gaps $g_n=(e_n^-,e_n^+)$}
\lb{z}
\end{figure}

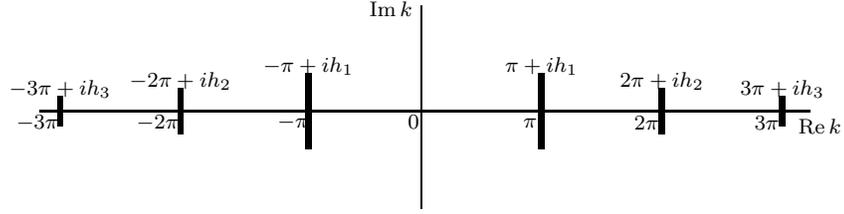
\begin{figure}
\tiny
\unitlength=1mm
\special{em:linewidth 0.4pt}
\linethickness{0.4pt}
\begin{picture}(120.67,34.33)
\put(20.33,20.00){\line(1,0){102.33}}
\put(71.00,7.00){\line(0,1){27.00}}
\put(70.00,18.67){\makebox(0,0)[cc]{$0$}}
\put(124.00,18.00){\makebox(0,0)[cc]{$\Re k$}}
\put(67.00,33.67){\makebox(0,0)[cc]{$\Im k$}}
\put(87.00,15.00){\linethickness{2.0pt}\line(0,1){10.}}
\put(103.00,17.00){\linethickness{2.0pt}\line(0,1){6.}}
\put(119.00,18.00){\linethickness{2.0pt}\line(0,1){4.}}
\put(56.00,15.00){\linethickness{2.0pt}\line(0,1){10.}}
\put(39.00,17.00){\linethickness{2.0pt}\line(0,1){6.}}
\put(23.00,18.00){\linethickness{2.0pt}\line(0,1){4.}}
\put(85.50,18.50){\makebox(0,0)[cc]{$\pi$}}
\put(54.00,18.50){\makebox(0,0)[cc]{$-\pi$}}
\put(101.00,18.50){\makebox(0,0)[cc]{$2\pi$}}
\put(36.00,18.50){\makebox(0,0)[cc]{$-2\pi$}}
\put(117.00,18.50){\makebox(0,0)[cc]{$3\pi$}}
\put(20.00,18.50){\makebox(0,0)[cc]{$-3\pi$}}
\put(87.00,26.00){\makebox(0,0)[cc]{$\pi+ih_1$}}
\put(56.00,26.00){\makebox(0,0)[cc]{$-\pi+ih_1$}}
\put(103.00,24.00){\makebox(0,0)[cc]{$2\pi+ih_2$}}
\put(39.00,24.00){\makebox(0,0)[cc]{$-2\pi+ih_2$}}
\put(119.00,23.00){\makebox(0,0)[cc]{$3\pi+ih_3$}}
\put(23.00,23.00){\makebox(0,0)[cc]{$-3\pi+ih_3$}}
\end{picture}
\caption{$k$-plane and cuts $\G_n=(\pi n-ih_n,\pi n+ih_n), \ n\in\Z$}
\lb{k}
\end{figure}

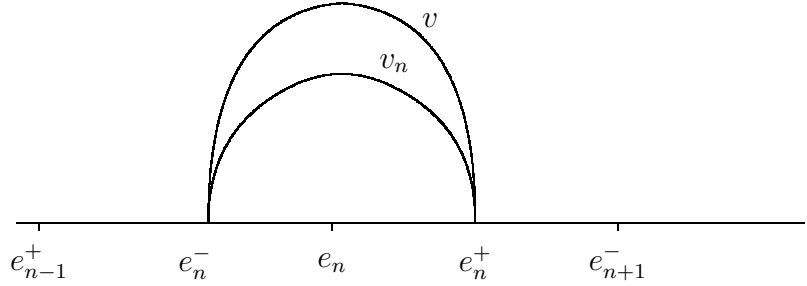
\begin{figure}
\unitlength 1mm 
\linethickness{0.4pt}
\ifx\plotpoint\undefined\newsavebox{\plotpoint}\fi 
\begin{picture}(119.75,66.5)(0,0)
\put(15,17.25){\line(1,0){104.75}}
\qbezier(40.5,17.25)(40.5,30.375)(52.25,35.75)
\qbezier(52.25,35.75)(58.25,38.375)(64.25,35.75)
\qbezier(64.25,35.75)(76.17,30.375)(76,17.25)
\qbezier(40.5,17.25)(40.5,45.125)(58.25,46.5)
\qbezier(76,17.25)(76.17,45.125)(58.25,46.5)
\put(38.75,12.00){\makebox(0,0)[cc]{$e_n^-$}}
\put(76,12){\makebox(0,0)[cc]{$e_n^+$}}
\put(65.25,38.5){\makebox(0,0)[cc]{$v_n$}}
\put(70,44.25){\makebox(0,0)[cc]{$v$}}
\put(18,17.25){\line(0,-1){1.00}}
\put(95,17.25){\line(0,-1){1.00}}
\put(57,17.25){\line(0,-1){1.00}}
\put(18,12.00){\makebox(0,0)[cc]{$e_{n-1}^+$}}
\put(95,12.00){\makebox(0,0)[cc]{$e_{n+1}^-$}}
\put(57,12.00){\makebox(0,0)[cc]{$e_{n}$}}
\end{picture}
\caption{The graph of $v(z+i0), \ z\in g_n\cup \s_n\cup \s_{n+1}$
and $h_n=v(e_n+i0)>0$}
\lb{grafv}
\end{figure}

Recall the identity from \cite{KK}. For each $n \in \Z$
the following identity holds true:
\[
\begin{aligned}
\lb{PIg}
&v(z+i0)=v_n(z)(1+Y_n(z)),\qqq \forall z\in g_n,
 \\
& \qq v_n(z)=|(z-e_n^+)(z-e_n^-)|^{1\/2},\qqq
Y_n(z)={1\/\pi}\int_{\R\sm g_n}{v(t)dt\/v_n(t)|t-z|}.
\end{aligned}
\]

\begin{lemma}
\lb{TY}
Let $Q_{2m}<\iy$ for some $m\ge 0$ and $s=\min_{n\ge 1}
|\s_n|$ and $M_n={1\/\pi}\int_{g_n}v(x)dx, n\in\Z$.
Then each function $Y_n, n\ge 1$, satisfies
\[
\lb{Y0}
Y_n^0:=\max_{z\in g_n} Y_n(z)\le  \sum_{j\ne n}
{M_j\/s^2|n-j|^2}\le {Q_0\/s^2},\qqq  if \qqq   m=0,
\]
\[
\lb{Y1}
\qqq \qqq \qqq\qqq \qqq \qqq   Y_n^0\le   {4Q_2\/n^2s^{4}},\qqq
  if \qqq     m\ge 1.
\]

\end{lemma}
\no {\bf Proof.} Using the estimate $\dist \{g_n, g_j\}\ge  s|n-j|$ we obtain
$$
Y_n(z)={1\/\pi}\int_{g\sm g_n}{v(t)dt\/v_n(t)|t-z|}=
\sum_{j\ne n} {1\/\pi}\int_{g_j}{v(t)dt\/v_n(t)|t-z|}
\le \sum_{j\ne n}  {1\/\pi}\int_{g_j}{v(t)dt\/s^2|n-j|^2}
=\sum_{j\ne n}  {M_j\/s^2|n-j|^2},
$$
which gives \er{Y0}. If $m\ge 1$, then the above estimates and
${1\/|j||n-j|}\le {2\/|n|}, j\ne n$ yield
$$
Y_n(z)\le \sum_{j\ne n}  {1\/\pi}\int_{g_j}{t^{2}v(t)dt\/s^2|n-j|^2t^{2}}\le
\sum_{j\ne n}  {1\/\pi}\int_{g_j}{t^{2}v(t)dt\/s^{4}|n-j|^2j^{2}}
\le
\sum_{j\ne n}  {4\/n^2\pi  s^{4}}\int_{g_j}t^{2}v(t)dt={4Q_2\/n^2s^{4}}.
$$
\BBox

We prove the main technical lemma of our paper.

\begin{lemma}
\lb{T2.1}
i)  Let $Q_{2m}<\iy$ for
some $m\ge 0$. Then the  quasimomentum has the form
\[
\lb{km1} k(z)=z-K_{{m}-1}(z)+f_m(z), \qqq \forall \ z\in \cZ,
\]
where
\[
\lb{km2}
f_m(z)={k_m(z)\/z^{2m}},\qqq
k_m(z)={1\/\pi}\int_{g}{t^{2m}v(t+i0)dt\/ t-z}.
\]
ii) Moreover, the following estimates and asymptotics hold true:
\[
\lb{km3}
|k_m(z)|\le {Q_{2m}\/\dist \{z,g\}},\qqq \forall z\in \cZ.
\]
\[
\lb{km4}
\max_{z\in g_n}|\Im f_m(z\pm i0)|\le h_n,  \qqq
\max_{z\in g_n}|\Re f_m(z\pm i0)|\le \max_{\pm} |f_m(e_n^\pm)|,
\]
\[
\lb{C1}
  f_m(e_n^\pm)= \Re f_m(e_n^\pm)=\mp{|g_n|\/2}(1+O(Y_n^0)),
\]
\[
\lb{C2}
\max_{z\in g_n} |f_m(z\pm i0)|=|g_n|(1+O(Y_n^0))
\]
as $n\to \iy$, uniformly in $z\in \ol g_n$, where
$Y_n^0=\max_{z\in g_n} Y_n(z)$.
\end{lemma}
\no {\bf Proof.} i) We have the simple identity
$$
{1\/ t-z}={1\/z^{2m}}{z^{2m}\/(t-z)}
={1\/z^{2m}}{z^{2m}-t^{2m}\/(t-z)}+ {1\/z^{2m}}{t^{2m}\/(t-z)}.
$$
Using this identity, we rewrite \er{Kq} in the form
(here and below $v(t)=v(t+i0), t\in \R$)
$$
k(z)-z={1\/\pi}\int_\R {v(t)\/t-z}dt=
{1\/\pi z^{2m}}\int_\R {(z^{2m}-t^{2m})v(t)\/t-z}dt+
{1\/\pi z^{2m}}\int_\R {t^{2m}v(t)\/t-z}dt.
$$
which gives \er{km1}, \er{km2}, since $Q_j=0$ for each odd $j$.

ii)  The identity \er{km2} gives \er{km3}.

 Using  $k=z-K_{m-1}+f_m$, we obtain
\[
\lb{esfm}
0\le \Im k(z+i0)= \Im f_m(z+i0)=v(z+i0)\le h_n, \qqq  \ \ z\in g_n.
\]
   Now we estimate the real part $\Re f_m(z+i0) , z\in g_n$.
Using $\Re k(z+i0)=\pi n$ on $g_n$, we obtain
\[
\begin{aligned}
\lb{fmder}
0=\Re  k'(z\pm i0)=1-K_{m-1}'(z)+\Re f_m(z\pm i0)'; \\
 \Re f_m(z\pm i0)'=-1+K_{m-1}'(z)<-1.
\end{aligned}
\]
Then the function ${\rm Re} f_m(x\pm i0)$ is decreasing in $x\in g_n$,
which yields  \er{km4}.

We prove \er{C1} for the case $e_n^-$. The proof for  $e_n^+$ is similar.
 Using \er{km2} we rewrite $f_m$ in the form
$$
f_m=f_{m1}+f_{m2},\qqq  f_{m1}(z)
={1\/\pi z^{2m}}\int_{g_n}{t^{2m}v(t)dt\/ t-z},\qqq
f_{m2}(z)={1\/\pi z^{2m}}\int_{g\sm g_n}{t^{2m}v(t)dt\/t-z},
$$
where $v(t)=v(t+i0)$.  Then using \er{PIg} and the new variable
$t=e_n^-+s$ we obtain
$$
f_{m1}(e_n^-)=
{1\/\pi }\int_0^{|g_n|}\rt(1+{s\/e_n^-}\rt)^{2m}{v_n(e_n^-+s)\/s}
(1+Y_n(e_n^-+s))ds
$$
$$
={1\/\pi}\int_0^{|g_n|}\rt(1+{O(s)\/e_n^-}
\rt)\rt|{|g_n|-s\/s}\rt|^{1\/2}(1+Y_n(t))ds =I_0+I_1,
$$
where
\[
I_0={1\/\pi}\int_0^{|g_n|}\rt(1+{O(s)\/e_n^-} \rt)\sqrt{|g_n|-s\/|s|}ds,\qqq
I_1={1\/\pi}\int_0^{|g_n|}\rt(1+{O(s)\/e_n^-} \rt)\sqrt{|g_n|-s\/|s|}Y_n(t)ds
\]
We have
\[
\lb{I01}
{1\/\pi}\int_0^{|g_n|}\sqrt{|g_n|-s\/s}ds
={|g_n|\/\pi}\int_0^{1}\sqrt{1-s\/s}ds={|g_n|\/2},
\]
and for the second term (in the case $m\ge 1$)  we have
$$
{1\/\pi e_n^-}\int_0^{|g_n|}{\sqrt{s(|g_n|-s)}}ds
={|g_n|^2\/\pi e_n^-}\int_0^{1}\sqrt{s(1-s)}ds
={|g_n|^2\/2 e_n^-}
$$
which yields
\[
\lb{I02}
I_0={|g_n|\/2}\rt(1+{O(|g_n|)\/n}\rt).
\]
Next, we consider $I_1$. Using \er{I01} we have
\[
I_1=  {1\/\pi}\int_0^{|g_n|}\sqrt{|g_n|-s\/s}dsO(Y_n^0 )
=|g_n|O(Y_n^0),\qqq   Y_n^0=\max_{t\in g_n}Y_n(t),
\]
which together with Lemma \ref{TY} yields \er{C1}.
In order to study $I_1$ we need to consider $Y_n$.
\BBox

{\bf Proof of Theorem \ref{T1}.}
Identities \er{km1} and \er{km2} imply  \er{13}, which yields \er{ask2}.
Thus we have
\[
\lb{13a}
k=z-K_{m}(z)+{k_{m+1}(z)\/z^{2m+2}},\qqq
k_{m+1}(z)={1\/\pi}\int_{\R}{t^{2m+2}v(t)dt\/ t-z},\qqq \ z\in \cZ.
\]
In order to show \er{ask1} we recall
the well known Nevanlinna Theorem (see \cite{Ah}).

\no   {\it Let $\m$ be a Borel measure on
$\R$ such that  $\int_{\R}(1+x^{2m})d\m(x)<+\iy$ for some
$m\ge 0$.  Then for each $A>0$ the following asymptotics hold true:}
$$
\begin{aligned}
&\int_{\R}{d\m(t)\/t-z}
=-\sum_{k=0}^{2m}{q_k\/z^{k+1}}+o({1\/z^{2p+1}})
\qq \as \qq|z|\to\iy, \qq y>A|x|,
 \end{aligned}
 $$
 where $q_j=\int_{\R}x^jd\m(x),\ 0\le
j\le 2m.$
Applying  Nevanlinna Theorem to \er{13a} and using \er{QmP} we obtain \er{ask1}.

The asymptotics \er{C2} and \er{km3}  give
$$
|f_{m+1}(z)|\le {|\g_n|\/2\pi n}+b_n,\qqq
b_n={O(|\g_n|)\/n^3}+{O(1)\/n^{2m+2}}\qqq z\in
\pa U_n,
$$
where $U_n=\{ z\in \cZ:  \dist \{z,g_n\} \le \ve\}.$
This yields \er{ask3} since the function $f_{m+1}$
is analytic in $U_n$.

The  asymptotics  \er{ask4} have been proved in \er{C1}.
$\BBox$.

\

\section {Asymptotics of the fundamental solutions}
\setcounter{equation}{0}

\

We recall some known facts from Lemma 3.1 in \cite{MO}. Define the
solution $y$ of the equation $-y''+qy=z^2y, z\ne 0$ in the form
\[
\begin{aligned}
\lb{y1}
y(x,z)=e^{\vk (x,z)},\ \ \ \ \ \vk (x,z)=iz x+\int_0^x\vk_* (t,z)dt,\\
\ \ \ \vk_*(x,z)=\sum_1^{m}{\vk_j(x)\/(2iz)^j}+{\vk_m(x,z)\/(2iz)^{m}},\\
y(0,z)=1,\qqq y'(0,z)=\vk'(0,z),
\end{aligned}
\]
for some $m\ge 1$. The function $\vk_*$ satisfies the equation:
\[
\lb{y2}
(2iz)\vk_*+\vk_*'+\vk_*^2=p.
\]
Moreover, the coefficients $\vk_{j}$ satisfy the following systems:
\[
\lb{y3}
 \vk_{j+1}=-\vk_{j}'-\sum_1^{j-1} \vk_{j-s}\vk_s,\ \ j\ge 1,
\]
where
\[
\begin{aligned}
\lb{y4}
\vk_1=p,\ \ \vk_2=-p',\ \ \   \vk_3=p''-p^2,\ \ \
\vk_4=-p'''+4pp',\ \ \ .... \\
\vk_j=(-1)^{j-1}p^{{j-1}}+\cP_{j-3}, \qq j=1,..., m-1,
\end{aligned}
\]
where $\cP_{j}$ is a polynomial in $p, p', p'', ..., p^{{j}}$.
The remainder $ \vk_{m}(x,z)$
satisfies
\[
\begin{aligned}
\lb{y5}
 \vk_{m}(0,z)= \vk_{m}(0,z)=0,\\
\vk_{m}(x,z)=O(1),\qqq \vk_{m}'(x,z)=O(z)
\end{aligned}
\]
as $|z|\to \iy$ uniformly in $[0,1]\ts \{z\in \C:\  |\Im z|\le r\}$
for any $r\ge 1$.
Define the functions
\[
\begin{aligned}
\lb{y6}
\r(z)\ev \vk'(0,z)=iz+\sum_1^{m}{\vk_j(0)\/(2iz)^j}\ev \t(z)+i\o(z),\\
\o(z)={\r(z)-\r(-z)\/2i},\ \ \  \t(z)={\r(z)+\r(-z)\/2},
\end{aligned}
\]
where
\[
\lb{ot}
\o(z)=z-\sum_0^{m-1\/2}{\vk_{2j+1}(0)\/(2z)^{2j+1}},\ \ \ \ \      \ \
\t(z)=\sum_1^{m\/2}{\vk_{2j}(0)\/(2z)^{2j}},
\]
\[
\lb{vk'1}
\vk'(1,z)=\r(z)+{\vk_m(1,z)\/(2iz)^m}.
\]
We rewrite the fundamental solutions $\vt, \vp$ in the forms
\[
\begin{aligned}
\lb{vpy}
\vp(x,z)={y(x,z)-y(x,-z)\/2i\o(z)},\\
\vp'(x,z)={y(x,z)\vk'(x,z)-y(x,-z)\vk'(x,-z)\/2i\o(z)},\
\end{aligned}
\]
and
\[
\lb{vty}
\begin{aligned}
\vt(x,z)={y(x,-z)\r(z)-y(x,z)\r(-z)\/2i\o(z)},\\
\vt'(x,z)={y(x,-z)\vk'(x,-z)\r(z)-y(x,z)\vk'(x,z)\r(-z)\/2i\o(z)}.
\end{aligned}
\]
Note that in \er{y1}-\er{vty} we do not use the condition $E_0^+=0$.
Recall that the set $\cZ(\ve)$ is given by
$\{z\in \cZ, \dist \{z,g\}> \ve\}, \ve >0$.

\begin{lemma}
\lb{T3.1}
 Let $p\in \mH_m$ for some $m\ge 0$ and let $r\ge 1$.
Then  the following asymptotics hold true:
\[
\lb{3as1}
\D(z)={y(1,z)+y(1,-z)\/2}+O(z^{-m})=\cos \x_m(1,z)+O(z^{-m}),
\]
\[
\lb{3as2}
\x(z):=\x_m(1,z)
=z-\sum_{0\le j\le {m-1\/2}}(-1)^{j}\int_0^1{\vk_{2j+1}(t)\/(2z)^{2j+1}}dt,
\]
\[
\lb{3as4}
y(1,z)=e^{i\x(z)+O(z^{-m})},
\]
and if in addition $E_0^+=0$, then
\[
\begin{aligned}
\lb{3as3}
k(z)=\x(z)+O(z^{-m}),\ \
\\
{(-1)^{j}\/2^{2j+1}}\int_0^1\vk_{2j+1}(t)dt=Q_{2j},\ \ \ \ \
\int_0^1\vk_{2j}(t)dt=0,
\ \ \  j\ge 0,
\end{aligned}
\]
as $|\Im z|\le r$ and $|z|\to \iy$.

\end{lemma}
{\bf Proof.} Let $A(z)=\vk'(1,z)-\r(-z)$. Identities \er{y6}-\er{ot}
yield
$$
A(z)=2i \o(z)+{\vk_m(1,z)\/(2iz)^m}=2i \o(z)+O(z^{-m}).
$$
Then this asymptotics and \er{vpy}, \er{vty} imply
$$
\D(z)={1\/4i\o(z)}
\lt(y(1,z)\vk'(1,z)-y(1,-z)\vk'(1,-z)+y(1,-z)\r(z)-y(1,z)\r(-z)\rt)
$$
$$
{y(1,z)A(z)-y(1,-z)A(-z)\/4i\o(z)}
={y(1,z)+y(1,-z)\/2}+O(z^{-m})=\cos i\vk (1,z)+O(z^{-m})
$$
The function $\D$ is real on the real line, which gives \er{3as2}, \er{3as4}.

If $E_0^+=0$, then the identity $\D(z)=\cos k(z), z\in \cZ$ and the asymptotics \er{3as1}, \er{3as2}
and the asymptotic estimate \er{km3} imply  \er{3as3}.
$\BBox$

Below, we need:

\begin{lemma}
\lb{T3.2}
 Let $p\in \mH_m$ for some $m\ge 0$ and let $r\ge 1$. Then the following asymptotics hold true:
\[
\begin{aligned}
\lb{vpy1}
\vp (1,z)={\sin \x(z)\/ \o(z)}+O(z^{-m-1}), \\
\vp'(1,z)=\cos \x(z)+{\t(z)\/ \o(z)}\sin \x(z)+O(z^{-m}),
\end{aligned}
\]
and
\[
\begin{aligned}
\lb{vty1}
\vt (1,z)=\cos \x(z)-{\t(z)\/ \o(z)}\sin \x(z)+O(z^{-m}),\\
\vt'(1,z)=-{\r(z)\r(-z)\/\o(z)}\sin \x(z)+O(z^{-m+1}),
\end{aligned}
\]
and
\[
\lb{by1}
\b(z)={\t(z)\/ \o(z)}\sin \x(z)+O(z^{-m}),
\]
as $|\Im z|\le r$ and $|z|\to \iy$.
\end{lemma}
{\bf Proof.}  Substituting the asymptotics \er{3as4} into the identity
\er{vpy} we have
$$
\vp (1,z)={y(1,z)-y(1,-z)\/2i\o (z) }=
{\sin \x(z)+O(z^{-m})\/ \o(z)}={\sin \x(z)\/ \o(z)}+O(z^{-m-1}),
$$
and using additionally \er{vk'1}, \er{y6} we have
$$
\vp'(1,z)={y(1,z)\vk'(1,z)-y(1,-z)\vk'(1,-z)\/2i\o(z)}=
{y_m(1,z)\r(z)-y_m(1,-z)\r(-z)+O(z^{-m})\/2i\o(z)},
$$
which yields \er{vpy1}.  The proof of the asymptotics in \er{vty1} is similar.

Using the asymptotics \er{vpy1}-\er{vty1}, we obtain
$$
\b(z)={\vp'(1,z)-\vt(1,z)\/2}={\t(z)\/ \o(z)}\sin \x(z)+O(z^{-m}),
$$
which yields \er{by1}.
$\BBox$

We need the following identities
\[
\begin{aligned}
\lb{fsb}
\P_{\pm}(0,z)=1, \qqq  \P_{\pm}'(0,z)=M_\pm(z), \\
  \P_{\pm}(1,z)=e^{\pm ik(z)}, \qqq
  \P_{\pm}'(1,z)=e^{\pm ik(z)}M_\pm(z),\qq \forall \ z\in \cZ.
  \end{aligned}
\]

{\bf Proof of Theorem \ref{T2}.}
 Using \er{vpy1}, \er{by1} and $k(z)=\x(z)+O(z^{-m})$ (see \er{3as3}),
  we have
$$
M_{\pm}(z)={\b(z){\pm}i\sin k(z)\/\vp (1,z)}=
{{\t(z)\/ \o(z)}\sin k(z)+O(z^{-m})\pm i\sin k(z)
\/ {\sin k(z)\/ \o(z)}+O(z^{-m})}
$$
$$
={\sin k(z) \r(\pm z)+O(z^{1-m})\/ \sin k(z)+O(z^{1-m})}.
$$
Moreover, if $z\in \cZ_\ve$, then
\[
\lb{mas}
M_{\pm}(z)=\r(\pm z)+O(z^{1-m}),
\]
which yields \er{aM}.
Using \er{vpy}, \er{vty}, \er{mas} and \er{y1}, we obtain
$$
\P_+(x,z)={1\/2i\o (z)}\lt[y(x,-z) \r(z)-y(x,z)
\r(-z) +M_+(z)(y(x,z)- y(x,-z))\rt]
=
$$
$$
={1\/2i\o (z)}\lt[y(x,z) (\r(z)-\r(-z)) +O(z^{1-m})(y(x,z)- y(x,-z))\rt]
$$$$
=y(x,z)+O(z^{-m})(y(x,z)- y(x,-z))=y(x,z)+O(z^{-m})=e^{i\x(x,\pm z)}+O(z^{-m}).
$$
The proof for $\P_-$ is similar. This yields \er{aP}.

The asymptotics \er{aqu} were proved in Lemma \ref{T3.1}.
$\BBox$

\

\section {Asymptotics for the distributions}
\setcounter{equation}{0}

\

In this Section we will determine the asymptotics of
the quasimomentum for the Schr\"odinger operator $H$  acting
in the Hilbert space $L^2(\R)$, given by
$$
Hy=-y''+(c+p')y.
$$
Here $p$ is a 1-periodic function belonging to
the real Hilbert space $\mH_*$ given by
$$
\mH_*=\rt\{p\in L^2(0,1): \int_0^1p(x)dx=0\rt\},
$$
and $c$ is a real constant.
Thus, $p'$ is a 1-periodic distribution, if $p' \in L^2(\T)$,
and then $H$ corresponds to the Hill operator with
$L^2$-potential. The situation considered in this paper, i.e. $p\in L^2(\T)$,
corresponds to a much more singular case.

We recall the results about the spectral properties of $H$ from \cite{K1}.
The spectrum of $H$ is purely absolutely continuous and
consists of intervals $\gS_n=[E_{n-1}^+,E_n^-]$. These intervals are
separated by the gaps $\g_n=(E_n^+,E_n^+)$ of length $|\g_n|\ge 0$.
If a gap $\g_n$ is degenerate, i.e. $|\g_n|=0$, then the corresponding segments
$\s_n,\s_{n+1}$ merge. We choose the constant $c$ in a way that $E_0^+=0$.
All these facts are similar to the case of smooth potentials.

We can not introduce the standard fundamental solutions for the operator $H$,
since the perturbation $p'$ is very strong. Thus we need another representation
of $H$. Define the unitary transformation $\mU:L^2(\R,\e^2 dx)\to L^2(\R,dx)$
as multiplication by $\e$. Thus $H$ is unitarily equivalent to
$$
H_1y=\mU^{-1}H\mU y=-{1\/\e^2}(\e^2y')'+(c-q^2)y=-y''-2py'+(c-p^2)y,
\ \ \  \e=e^{\int_0^xp(t)dt}
$$
acting in $L^2(\R,\e^2 dx)$. This representation is more convenient,
since we can introduce the fundamental solutions
$\vp_1(x,z),\vt_1(x,z)$ of the equation
\[
\lb{1.1}
-y''-2qy'+(c-p^2)y=z^2 y, \ \ \ \ \ \ \ \ z\in\C,
\]
with the conditions: $\vp_1(0,z)=\vt_1'(0,z)=0,
\vp_1'(0,z)=\vt_1(0,z)=1$. Define the Lyapunov function
$$
\D(z)={\vp_1'(1,z)+\vt_1(1,z)\/2}.
$$
Similar to the case of smooth potentials, we introduce the
quasimomentum  $k(z)$ and the momentum domain $\cZ$ and the
quasimomentum domain $\cK$ by \er{5} and \er{cK}. The quasimomentum
$k(z)$ is a conformal mapping from $\cZ$ onto $\cK$, and it
satisfies the standard properties \er{pk} and \er{2.3}-\er{PIg}. To
characterize the quasimomentum $k(z)$ further, we recall the results
from \cite{K1}:
 The  quasimomentum has the form
\[
\lb{d1} k(z)=z-k_0(z), \qqq k_0(z)={1\/\pi}\int_{g}{v(t)dt\/
t-z},\qqq \forall \ z\in \cZ,
\]
and for any $A>0$, the following asymptotics hold true
\[
\lb{da}
k(z)=z-{P_{-1}+o(1)\/z}\qqq as \ |z|\to \iy, \qq y>A|x|.
\]
Here, the coefficient $P_{-1}$ has the form
\[
P_{-1}={\|q\|^2\/2}={1\/\pi}\int_\R v(t+i0)dt=
{1\/2\pi}\int\!\!\int_\C |k'(z)-1|^2dxdy, \ \ \ \ z=x+iy,
\]
where $q\in \mH_{*}$ is a solution of the Riccati  equation
\[
\lb{RE}
p'=q'(x)+q(x)^2-\|q\|^2.
\]
Recall that the mapping $p\to q$ acting from $\mH_*$ into
$\mH_*$ is a real analytic isomorphism onto itself.
Thus for each $p\in \mH_{*}$ there exists a unique solution
$q\in \mH_*$ of the equation \er{RE}.

\medskip

Define the sequence $S_n, n\ge 1$ by
\[
S_n(r)=\sum_{j\in \Z}{M_j\/|n-j|_1}, \qqq |j|_1=\ca s|j|, &if \ j\ne 0\\
{r\/2},  &if\   j= 0 \ac,\qq M_j={1\/\pi}\int_{g_j}v(t+i0)dt.
\]
Note that simple estimates yield
$$
\sum_{n\ge 1} S_n^a(r)\le
\sum_{n\ge 1}\sum_{j\in \Z}{M_{j}\/|n-j|_1^{a}}
\rt(\sum_{j'\in \Z}M_{j'}\rt)^{a-1}
\le Q_0^a\rt({4^a\/r^{2a}}+{2\/s^a}C_a\rt),
$$
where $C_a=\sum_{n\ge 1}{1\/|j|^{a}}$. This yields
\[
\lb{esSn}
\sum_{n\ge 1} S_n^a(r)\le 4Q_0^a\rt({1\/r^2}+{1\/s^2}\rt), \qqq if \qq a>1.
\]

Define the domains
\[
\lb{Vn}
V_n(r)=\rt\{ z\in \C:  |\Im  z|\le r, \qq
{e_n^-+e_{n-1}^+\/2}< \Re z < {e_n^++e_{n+1}^-\/2}\rt\}\sm U_n, \ r>\pi,
\]
\[
\lb{Un}
U_n=\rt\{ z\in \cZ:  \dist \{z,g_n\} \le \ve\rt\}.
\]

\begin{theorem}   \lb{Td}
 Let $H=-{d^2\/dx^2}+(c+p')$ for some $p\in \mH_0$,
 and let $r\ge \pi$. Then for each $\ve >0$ small enough,
 $k_0(z)$ satisfies
\[
\lb{Q01}
\max_{z\in g_n} |k_0(z\pm i0)|=|g_n|(1+O(Y_n^0))
\]
\[
\lb{Q02}
\max_{\dist\{z,g_n\}=\ve}|k_0(z)|\le  S_n(\ve),  \qqq
\]
\[
\lb{Q03}
\max_{z\in \pa V_n\sm \pa U_n}|k_0(z)|\le  S_n(1),
\]
\[
\lb{Q04}
\max_{z\in U_n} |k_0(z\pm i0)|=|g_n|(1+O(Y_n^0)),
\]
\[
\lb{Q05}
\max_{z\in V_n}|k_0(z)|\le  S_n(\ve)+S_n(s),
\]
\end{theorem}
\no {\bf Proof.} The asymptotics \er{Q01} follows from \er{C2}.

In order to show   \er{Q02}, we write the simple estimate of $k_0$
in the form:
$$
|k_0(z)|\le {1\/\pi}\int_{g}{v(t)dt\/|t-z|}=
\sum_{j\in \Z}{1\/\pi}\int_{g_j}{v(t)dt\/|t-z|}.
$$
Consider the summands separately in sum. We obtain estimates for $z\in U_n$:
$$
{1\/\pi}\int_{g_n}{v(t)dt\/|t-z|}\le {1\/\pi \ve}\int_{g_n}v(t)dt=
 {M_n\/\ve}, \qqq n=j;
$$
$$
{1\/\pi}\int_{g_j}{v(t)dt\/|t-z|}\le {1\/\pi s}\int_{g_j}{v(t)dt\/|n-j|}
= {M_j\/s|n-j|}, \qqq n\ne j.
$$
Summing these estimates we obtain \er{Q02}. The proof of \er{Q03} is similar.

 The function $k_0$ is analytic in the domain $U_n$ and satisfies
 the estimates \er{Q01}, \er{Q02} on the boundary, which yields \er{Q04}.

 The function $k_0$ is analytic in the domain $V_n$ and
 satisfies the estimates \er{Q02}, \er{Q03} on the boundary,
 which yields \er{Q05}.
\BBox


\no  {\bf Acknowledgments.}\small 
This work was supported by the
Ministry of education and science of the Russian Federation, state
contract 14.740.11.0581.

\end{document}